\documentclass[11pt]{article}
\usepackage{amssymb,amsfonts,amsmath,latexsym,epsf,tikz,url}

\newtheorem{theorem}{Theorem}[section]

\newtheorem{conjecture}[theorem]{Conjecture}
\newtheorem{corollary}[theorem]{Corollary}
\newtheorem{lemma}[theorem]{Lemma}

\newtheorem{example}[theorem]{Example}
\newcommand{\proof}{\noindent{\bf Proof. }}
\newcommand{\qed}{\hfill $\square$\medskip}

\begin{document}

\title{The distinguishing number and the distinguishing index of  Cayley graphs}

\author{
Saeid Alikhani  $^{}$\footnote{Corresponding author}
\and
Samaneh Soltani
}

\date{\today}

\maketitle

\begin{center}
Department of Mathematics, Yazd University, 89195-741, Yazd, Iran\\
{\tt alikhani@yazd.ac.ir, s.soltani1979@gmail.com}
\end{center}


\begin{abstract}
The distinguishing number (index) $D(G)$ ($D'(G)$) of a graph $G$ is the least integer $d$
such that $G$ has an vertex labeling (edge labeling)  with $d$ labels  that is preserved only by a trivial
automorphism. In this paper, we investigate the distinguishing number and the distinguishing  index of Cayley graphs.
\end{abstract}


\medskip
\noindent{\bf AMS Subj.\ Class.}: 05C25, 05E18  

\section{Introduction and definitions}
Let $\Gamma = (V,E)$ be a simple and undirected graph. The set of all automorphisms of $\Gamma$ forms a permutation group called the \textit{automorphism group}  of $\Gamma$, which we denote by
${\rm Aut}(\Gamma)$. Given a group $H$ and a subset $S \subseteq H$, the \textit{Cayley digraph}  of $H$ with respect
to $S$, denoted by ${\rm Cay}(H, S)$, is the digraph with vertex set $H$ and with an arc from $h$
to $sh$ whenever $h \in H$ and $s \in  S$. When $S$ is closed under inverses, $(g, h)$ is an arc of
the Cayley digraph if and only if $(h, g)$ is an arc, and so we can identify the two arcs
$(g, h)$ and $(h, g)$ with the undirected edge $\{g, h\}$. When $1 \notin S$, the Cayley digraph
contains no self-loops. Thus, when $1 \notin  S = S^{-1}$, ${\rm Cay}(H, S)$ can be considered to be
a simple, undirected graph. The Cayley graph ${\rm Cay}(H, S)$ is connected if and only if
$S$ generates $H$. 

The automorphism group of the Cayley graph ${\rm Cay}(H, S)$ contains the right regular
representation $R(H)$ as a subgroup, and hence all Cayley graphs are vertex-transitive.
Let $e$ denote the identity element of the group $H$ and also the corresponding vertex
of ${\rm Cay}(H, S)$. Since $R(H)$ is regular, ${\rm Aut}(\Gamma) = {\rm Aut}(\Gamma)_eR(H)$, where ${\rm Aut}(\Gamma)_e$ is the
stabilizer of e in ${\rm Aut}(\Gamma)$. The set of automorphisms of the group $H$ that fixes $S$
setwise, denoted by ${\rm Aut}(H, S) := \{\Gamma \in  {\rm Aut}(H) : S^{\Gamma} = S\}$, is a subgroup of ${\rm Aut}(\Gamma)_e$ (see \cite{N. L. Biggs}). Godsil in \cite{C. Godsil}, and Xu in \cite{M. Y. Xu} showed that for any Cayley graph $\Gamma = {\rm Cay}(H, S)$, the normalizer $N_{{\rm Aut}}(\Gamma)(R(H))$ is
equal to the semidirect product $R(H) \rtimes {\rm Aut}(H, S)$. A Cayley graph
$\Gamma := {\rm Cay}(H, S)$ is said to be \textit{normal} if $R(H)$ is a normal subgroup of ${\rm Aut}(\Gamma)$, or
equivalently, if ${\rm Aut}(\Gamma) = R(H) \rtimes {\rm Aut}(H, S)$. Thus, normal Cayley graphs are those  have the smallest possible full automorphism group (\cite{M. Y. Xu}).

Let $G$ be a simple graph. A labeling of $G$, $\phi : V \rightarrow \{1, 2, \ldots , r\}$, is said to be \textit{$r$-distinguishing},  if no non-trivial  automorphism of $G$ preserves all of the vertex labels. The point of the labels on the vertices is to destroy the symmetries of the
graph, that is, to make the automorphism group of the labeled graph trivial.
Formally, $\phi$ is $r$-distinguishing if for every non-identity $\sigma \in {\rm Aut}(G)$, there exists $x$ in $V$ such that $\phi(x) \neq \phi(\sigma(x))$. The \textit{distinguishing number} of a graph $G$ is defined  by
\begin{equation*}
D(G) = {\rm min}\{r \vert ~ G ~\text{{\rm has a labeling that is $r$-distinguishing}}\}.
\end{equation*} 

This number has defined by Albertson and Collins \cite{Albert}. Similar to this definition, Kalinowski and Pil\'sniak \cite{Kali1} have defined the \textit{distinguishing  index} $D'(G)$ of a graph $G$ which is  the least integer $d$ such that $G$ has an
edge labeling with $d$ labels that is preserved only by the identity automorphism of $G$. If a graph has no nontrivial automorphisms, its distinguishing number is one. In other words, $D(G) = 1$ for the asymmetric graphs.
 The other extreme, $D(G) = \vert V(G) \vert$, occurs if and only if $G = K_n$. The distinguishing index of some examples of graphs was exhibited in \cite{Kali1}. In the sequel, we need the following results:
 
\begin{theorem}{\rm \cite{K.L. Collins A.N. Trenk,klavzaretal}}\label{disnumbound}  If $G$ is a connected graph with maximum degree $\Delta$, then $D(G) \leq \Delta+1$. Furthermore,
equality holds if and only if $G$ is a $K_n$, $K_{n,n}$, $C_3$, $C_4$ or $C_5$.
 \end{theorem}
 \begin{theorem}{\rm \cite{pilsniaknord}}\label{distindex}
 Let $G$ be a connected graph that is neither a symmetric nor an asymmetric tree. If the maximum degree of $G$ is at least $3$, then $D'(G) \leq \Delta(G) - 1$ unless $G$ is $K_4$ or $K_{3,3}$.
  \end{theorem}
 
 The concept was naturally extended to general group actions \cite{Tymockzko}. 
Let $\Gamma$ be a group acting on a set $X$. If $g$ is an element of $\Gamma$ and $x$ is in $X$ then denote the action of $g$ on $x$ by $g.x$. Write $\Gamma .x$ for the orbit containing $x$. Recall that stabilizer of the subset $Y \subseteq X$ is defined to be ${\rm Stab}_{\Gamma}(Y ) = \{g \in \Gamma : g.y = y ~{\rm for~ all}~ y \in Y \}$. We sometimes omit the subscript and write ${\rm Stab}(Y)$. For a positive integer $r$, an \textit{$r$-labeling} of $X$ is an onto function $\phi : X\rightarrow \{1,2,\ldots ,r\}$. We say $\phi$ is a \textit{distinguishing labeling} (with respect to the action of $\Gamma$) if  the only group elements that preserve the labeling are in ${\rm Stab}_{\Gamma}(X)$. The \textit{distinguishing number} $D_{\Gamma}(X)$ of the action of  $\Gamma$ on $X$ is defined as 
\begin{equation*}
D_{\Gamma}(X) ={\rm min}\{r:~ {\rm there~ exists~ an~ r-distinguishing~ labeling}\}.
\end{equation*}
In particular, if the action is faithful, then the only element of $\Gamma$ preserving labels is the identity. 
 
 \section{The distinguishing number of Cayley graphs}
 
 We start this section with obtaining  general bounds for the distinguishing number of Cayley graphs.
 
 \begin{theorem}\label{up1}
 If $\Gamma = {\rm Cay}(H, S)$ is a Cayley graph, then $2 \leq D(\Gamma) \leq |S|+1$.
 \end{theorem}
\proof
Since $R(H) \subseteq {\rm Aut}(\Gamma)$, so $D(\Gamma)\geq 2$. On the other hand, $\Gamma$ is a regular graph of degree $|S|$, and hence $D(\Gamma) \leq |S|+1$, by Theorem \ref{disnumbound}.\qed
 
 If $\Gamma = {\rm Cay}(H, S)$ is a normal Cayley graph, then we can improve the upper bound in Theorem \ref{up1}.

 \begin{theorem}
 If $\Gamma = {\rm Cay}(H, S)$ is normal, then  $D(\Gamma) \leq D_{{\rm Aut}(H,S)}(S)+1$.
 \end{theorem}
\proof 
Since $\Gamma$ is normal, so ${\rm Aut}(\Gamma) = R(H) \rtimes {\rm Aut}(H, S)$. Set $D_{{\rm Aut}(H,S)}(S)=t$. We label the elements of $S$ distinguishingly with $t$ labels $\{1,2,\ldots , t\}$, next we label all elements of $G\setminus S$ with a new label, say $0$. This vertex labeling of $\Gamma$ is distinguishing, because if $f$ is an automorphism of $\Gamma$ preserving the labeling, then $f$ must map $G\setminus S$ to $G\setminus S$ and $S$ to $S$, setwise. Hence $f\in {\rm Aut}(H,S)$. Since the elements of $S$ have been labeled distinguishingly, so $f$ is the identity automorphism.\qed
 
 \begin{corollary}
  If $\Gamma = {\rm Cay}(H, S)$ is a normal Cayley graph such that   ${\rm Aut}(H,S)=\{id\}$, then the distinguishing number of $\Gamma$ is $D(\Gamma)=2$.
 \end{corollary}

 Let $G$ be an (abstract) finite group. If $G\cong {\rm Aut}(\Gamma)$ for some Cayley graph $\Gamma=Cay(H,S)$, then $H$ has to be isomorphic to a subgroup of $G$. To simplify our notation,  simply assume that $H$ is a subgroup of $G$. In the case when $H = G$, the action of $G$ on the vertices of $\Gamma$ is regular, and the Cayley graph $Cay(G,S)$ is called a \textit{graphical regular representation} of $G$. In other words, a Cayley graph  $\Gamma = {\rm Cay}(H, S)$  is called a graphical regular representation of the group G if ${\rm Aut}(\Gamma) = H$. Conder and Tucker in \cite{M. Conderetal} showed that if $A$ is a group acting on the set $X$, then  $D_A(X) = 2$ if and only if $A$ has a regular orbit on the subsets of $X$.  
 
 \begin{theorem}\label{graphregrep}
  If $\Gamma = {\rm Cay}(H, S)$ is a graphical regular representation of $H$, then the distinguishing number of $\Gamma$ is  $D(\Gamma) = 2$.
 \end{theorem}
\proof Since $\Gamma$ is a graphical regular representation of $H$, so ${\rm Aut} (\Gamma) = H$, and so the action of $G$ on the vertices of $\Gamma$ is regular. Thus $D(\Gamma) = 2$.  \qed

\begin{example}
As usual, we denote by ${\rm Sym}(\Omega)$ the group of all permutations
of a set $\Omega$, which is called the symmetric group
on $\Omega$. If $\Omega$ is the set $\{1, 2, \ldots  ,n\}$, then the symmetric group on $\Omega$  is also denoted by $S_n$. The {\rm Pancake graph} $P_n$,
also called the prefix-reversal graph, is the Cayley graph
${\rm Cay}(S_n, PR_n)$, where $PR_n =\{r_{1 j}~: 2\leq  j  \leq n\}$ and $r_{1j}$ is the following permutation:
\begin{equation*}
r_{1j}=\left(
\begin{array}{ccccccc}
1 & 2 & \ldots & j & j+1  & \ldots & n\\
j & j-1 & \ldots & 1 & j+1  & \ldots & n
\end{array}\right).
\end{equation*}

Deng and Zhang in \cite{Dengetal} showed that $P_n$ $(n \geq 5)$ is a graphical regular representation of $S_n$. Therefore by Theorem \ref{graphregrep}, for any $n \geq 5$,  $D(P_n) = 2$.  
\end{example}

Now we want to show that the distinguishing number of normal Cayley graphs is at most three with finitely many exceptions. For this purpose, we need the following theorem.

\begin{theorem}{\rm \cite{M. Conderetal}}\label{disgroup}
 For groups, $D_{{\rm Aut}(H)}(H) \neq 2$ if and only if $H$ is isomorphic to one of the
four elementary abelian groups $C_2^2$, $C_2^3$, $C_3^2$ and $C_2^4$ (of orders $4$, $8$, $9$ and $16$) or the quaternion group $Q_8$.
\end{theorem}

\begin{theorem}
If  $\Gamma = {\rm Cay}(H, S)$ is a normal Cayley graph such that  $H$ is not isomorphic to one of the
four elementary abelian groups $C_2^2$, $C_2^3$, $C_3^2$ and $C_2^4$ (of orders 4, 8, 9 and 16) or the quaternion group $Q_8$, then $D(\Gamma) \leq 3$.
\end{theorem}
\proof First by Theorem \ref{disgroup}, we label the vertex set of $\Gamma$ (elements of $H$) with the two labels $1$ and $2$ distinguishingly under ${\rm Aut}(H, S)$. Next we change the label of the vertex $e$, i.e., the identity element of $H$, to the new label $3$. To show this labeling is distinguishing, we use  normality of $\Gamma$. Since $\Gamma$ is normal, so ${\rm Aut}(\Gamma) = R(H) \rtimes {\rm Aut}(H, S)$. Hence if $f\in {\rm Aut}(\Gamma)$, then there exists $a \in H$ and $g \in {\rm Aut}(H,S)$ such that $f(x) = g(x).a$. If $f$ preserves the labeling then since $e$ is the only vertex of $\Gamma$ with label $3$, so $f(e) = e$. Now from $g(e)= e$, and also $f(e) = e$ we obtain that $a = e$. Thus any automorphism of $\Gamma$ preserving the labeling is an elements of ${\rm Aut}(H, S)$. Since we labeled the elements of $H$ such that the only automorphism of $H$ preserving the labeling is the identity, so $f$ is the identity automorphism. \qed

Let $S$ be a set of transpositions generating the symmetric group $S_n$. The \textit{transposition
graph}  of $S$, denoted by $T(S)$, is defined to be the graph with vertex set
$\{1, \ldots , n\}$, and with two vertices $i$ and $j$ being adjacent in $T(S)$ whenever $(i, j) \in  S$.
A set $S$ of transpositions generates $S_n$ if and only if the transposition graph of $S$ is
connected. The \textit{bubble-sort graph} of dimension $n$,  $B_n$, is the Cayley graph of $S_n$ with respect to the generator set $\{(1, 2), (2, 3), \ldots  , (n-1, n)\}$. In other words, the bubble-sort graph is the Cayley graph ${\rm Cay}(S_n, S)$ corresponding to the case where the transposition graph $T(S)$ is the path graph on $n$ vertices (\cite{A. Ganesan}).

\begin{theorem}{\rm \cite{A. Ganesan}}\label{autbubb}
	 Let $S$ be a set of transpositions generating $S_n$ $(n \geq 3)$ such that the Cayley graph ${\rm Cay}(S_n, S)$ is normal. Then the automorphism group of the Cayley graph
${\rm Cay}(S_n, S)$ is the direct product $S_n \times  {\rm Aut}(T(S))$, where $T(S)$ denotes the transposition graph of $S$.
\end{theorem}

In the following result we want to obtain the distinguishing number of bubble-sort graphs. Let $S$ be a set of transpositions generating $S_n$. Let $L$ denote the left regular action of $S_n$ on itself, defined by $L : S_n \rightarrow {\rm Sym}(S_n)$, $a \mapsto L(a)$, where $L(a) : x \mapsto a^{-1}x$.  The following lemma follows directly from the proof of Theorem \ref{autbubb}:

\begin{lemma}\label{newlemma}
	Every element in automorphism group of bubble-sort graph,   ${\rm Aut}(B_n)$, can be expressed uniquely in the form $R(a)L(b)$ for some $a\in S_n$ and $b\in {\rm Aut}(T(S))$.
\end{lemma} 

\begin{theorem}
The distinguishing number of bubble-sort graph $B_n$ of dimension $n\geq 3$ is two.
\end{theorem}
\proof   By Lemma \ref{newlemma}  every element in ${\rm Aut}(B_n)$ can be expressed uniquely in the form $R(a)L(b)$ for some $a\in S_n$ and $b\in {\rm Aut}(T(S))$. Since for the bubble-sort graphs, the graph $T(S)$ is a path graph of order $n$, so we denote the automorphism group of $T(S)$ as ${\rm Aut}(T(S))=\{id, \sigma\}$ where $\sigma= (1~n)(2~n-1)(3~n-2)(4~n-3) \cdots$. If $n=3$, then $B_3$ is a cycle graph of order six, and so $D(B_3) = 2$. Hence we suppose that $n \geq 4$. We present a labeling of vertices of $B_n$ (the elements of $S_n$) and then we show that  this labeling is distinguishing. We label the vertices $(1~2), (1~2~3), \ldots , (1~2~3~\cdots~n)$ with label $1$, and the remaining vertices of $B_n$ with label $2$. By  contradiction suppose that $f$ is a non-identity automorphism of $B_n$ preserving the labeling. Thus $f(x)=b^{-1} x a$ for some  $a\in S_n$ and $b\in \{id, \sigma\}$. We consider the two following cases:

\medskip
Case 1) Let $b= id$. Since $f$ preserves the labeling, then $f$ maps $(1~2)$ to one of elements of the set $\{(1~2), (1~2~3), \ldots , (1~2~3~\cdots~n)\}$. If $f((1~2))= (1~2)$, then $(1~2) a = (1~2)$, and hence $a = 1$. Therefore $f$ is the identity automorphism, which is a contradiction.  If $f((1~2))= (1~2~3~\cdots ~i)$ for any $i \geq 3$, then $(1~2) a = (1~2~3~\cdots ~i)$, and hence $a = (1~3~4~5~\cdots ~ i)$. In this case $f$ maps the vertex $(1~2~3~4)$ to the vertex $(1~2~4~3~5~6~\cdots ~i)$. Since the label of vertices $(1~2~3~4)$ and $(1~2~4~3~5~6~\cdots ~i)$ are different, so  $f$ dose not preserve the labeling, which is a contradiction. 
 
Case 2) Let $b = \sigma$. Since $f$ preserves the labeling, then $f$ maps $(1~2)$ to one of elements of the set $\{(1~2), (1~2~3), \ldots , (1~2~3~\cdots~n)\}$. If $f((1~2))= (1~2~\cdots ~i)$ for any $i\geq 2$, then it can be seen that $a$ is a permutation of $S_n$ such that $a(n) = 2$. Hence  $f$ maps the vertex $(1~2~3~\cdots ~n)$ to the vertex $x$ in $S_n$ such that $x(2) =2$, and so the label of $x$ is 2. Since the label of vertices $(1~2~3~\cdots ~n)$ and $x$ are different, so  $f$ dose not preserve the labeling, which is a contradiction.  \qed

 Now we state the   following theorem.
\begin{theorem}{\rm \cite{A.V. Konygin}}\label{primitive}
Let $\Gamma$ be a finite connected undirected graph without loops or multiple edges,
and let $\Gamma$ admit a vertex-primitive group of automorphisms. Then either $D(\Gamma) = 2$ or one of the
following assertions holds:
\begin{enumerate}
\item[(i)] $\Gamma$ is a complete graph and $D(\Gamma) = |V (\Gamma)|$;
\item[(ii)] $D(\Gamma) = 3$ and $\Gamma$ is isomorphic to one of the following four graphs: the cycle of length 5; the
Petersen graph; the graph complementary to the Petersen graph; the graph with the set of vertices
$\{(i, j) ~|~ i, j \in  \{1, 2, 3\}\}$ in which the vertex $(i, j)$ is adjacent to the vertex $(i', j')$ if $i = i'$ or $j = j'$.
\end{enumerate}
\end{theorem}

It is clear that the automorphism group of the   Cayley $\Gamma = {\rm Cay}(H, S)$ acts transitivity on the vertex set $H$. If this action is primitive then we can obtain the distinguishing number of such graphs as follows (directly by  Theorem \ref{primitive}).
\begin{theorem}\label{thmprimimprim}
 Let  $\Gamma = {\rm Cay}(H, S)$ be a Cayley  graph.
  If ${\rm Aut}(\Gamma)$ acts on $H$ primitively and $\Gamma$ is not a complete graph, then $D(\Gamma)=2$. If  $\Gamma = {\rm Cay}(H, S)$ is a complete graph then $D(\Gamma) = |H|+1$.

\end{theorem}

 \section{The distinguishing index of Cayley graphs}
 
 In this section we study the distinguishing index of Cayley graphs.  
 In 1969, Lov\'asz \cite{lovasz}  asked whether every finite connected vertex-transitive graph has a Hamilton path, that is, a simple path going through all vertices. The graphs with a Hamiltonian path are called \textit{traceable graphs}. Also,    Parsons et.al. conjectured that there is a Hamiltonian cycle in every
Cayley graph. If this conjecture is true, then by the following theorem, we can conclude  that the distinguishing index of all Cayley graphs of order at least seven, is at most two. 
 
 \begin{theorem}{\rm \cite{pilsniaknord}}\label{traceble}
 If $G$ is a traceable graph of order $n\geq  7$, then $D'(G) \leq 2$.
 \end{theorem}
 
 \begin{theorem}\label{restricindexcayl}
 Let $\Gamma$ be a $k$-regular graph of order $n$.
 \begin{itemize}
 \item[(i)] If $k \leq 4$, then $D'(\Gamma) \leq 3$.
 \item[(ii)] If $k\geq \frac{n-1}{2} \geq 3$, then $D'(\Gamma) =2$.
\item[(iii)] Let $\Gamma = {\rm Cay}(H, S)$ be a Cayley  graph of order $|H|\geq 2$. If   $\Gamma = {\rm Cay}(H, S)$ is a graphical regular representation of $H$, then $D'(\Gamma) =2$.
 \end{itemize}
 \end{theorem}
 \proof 
  \begin{itemize}
  	\item[(i)] The result follows from Theorem \ref{distindex}.
  	\item[(ii)] Since $k\geq \frac{n-1}{2}$, so $\Gamma$ is a traceable graph and  since $k \geq 3$, so the order of $\Gamma$ is greater or equal than seven, and hence  the result follows from Theorem \ref{traceble}. 
  	
  	\item[(iii)] Since $\Gamma$ is a  graphical regular representation of $H$, so ${\rm Aut}(\Gamma) = R(H)$.  Let $S=\{s_1,s_2, \ldots , s_{|S|}\}$. We label the edges $(e,s_i)$ with label $1$ for any $1\leq i \leq |S|$, and the remaining edges with label $2$ ($e$ is the identity element of the group $H$). Let $f\in R(H)$ be an automorphism of $\Gamma$. Thus there exists $a\in H$ such that $f(x) = x.a$ for any $x\in H$.  If $f$ preserves the labeling, then since $e$ is the only vertex such that all its incident edges have the label $1$, so  $f$ must map the vertex $e$ to $e$. Hence $a=e$, and therefore $f$ is the identity automorphism. Thus this $2$-labeling is distinguishing. Since $|{\rm Aut}(\Gamma)| > 1$, so $D'(\Gamma) = 2$.\qed 
  \end{itemize}

 In the following result which is one of the main result of this paper,   we show that the distinguishing index of regular graphs is at most three.
 \begin{theorem}
 If $G$ is a connected $k$-regular graph of order $n$, then $D'(G)\leq 3$.
 \end{theorem}
 \proof   If $k \leq 4$ or $k\geq \frac{n-1}{2} \geq 3$,  then the result follows from Theorem \ref{restricindexcayl}. Then, we suppose that $5\leq k <  \frac{n-1}{2}$. 
   Let $v$ be a vertex of $G$ with the maximum degree $\Delta$.   Let $N^{(1)}(v)= N_G(v) =\{v_1, \ldots , v_k\}$ be the vertices of $G$ at distance one from $v$.  By the following steps, we label  the edges of graph: 

Step 1) We can label the edges $vv_{i (\lceil \sqrt[k]{k^{k -1}}\rceil -1) +j}$ with label $i$, for $0 \leq i \leq 2$ and $1 \leq j \leq \lceil \sqrt[k]{k^{k -1}}\rceil -1$, and we do not use label 0 any more. With respect to the number of incident edges to $v$ with label 0, we conclude that the vertex $v$ is fixed under each automorphism of $G$ preserving the labeling.  Hence, every automorphism of $G$ preserving the labeling must map the set of vertices of $G$ at distance $i$ from $v$ to itself setwise,  for any $1 \leq i \leq {\rm diam}(G)$. We denote the set of vertices of $G$ at distance $i$ from $v$  for any $2 \leq i \leq {\rm diam}(G)$, by $N^{(i)}(v)$. If $N^{(i)}(v) = \emptyset$, for any $i \geq 2$, then $G$ has a Hamiltonian path, and since $k \geq 5$, so the order of $G$ is at least 7, and hence $D'(G) \leq 2$ by Theorem \ref{traceble}. Thus we suppose that $N^{(i)}(v) \neq  \emptyset$, for some $i \geq 2$. 

 Now we partition the vertices $N^{(1)}(v)$ to two sets $M_1^{(1)}$ and $M_2^{(1)}$ as follows:
 \begin{equation*}
 M_1^{(1)} =\left\{x\in N^{(1)}(v)~:~ N(x) \subseteq N(v)\right\},~~M_2^{(1)}=\left\{x\in N^{(1)}(v)~:~ N(x) \nsubseteq N(v)\right\}.
 \end{equation*}
 Thus the sets $M_1^{(1)}$ and $M_2^{(1)}$ are mapped to $M_1^{(1)}$ and $M_2^{(1)}$, respectively, setwise, under each automorphism of $G$ preserving the labeling. For  $0 \leq i \leq 2$, we set $L_i =\{v_{i (\lceil \sqrt[k]{k^{k -1}}\rceil -1)+j}~:~ 1 \leq j \leq \lceil \sqrt[k]{\Delta^{k -1}}\rceil-1\}$. By this notation, we get that for $0 \leq i \leq 2$, the set $L_i$ is mapped to $L_i$  under each automorphism of $G$ preserving the labeling, setwise. Let the sets $M_{1i}^{(1)}$ and $M_{2i}^{(1)}$ for $0 \leq i \leq 2$ are as follows:
\begin{equation*}
M_{1i}^{(1)}= M_1^{(1)} \cap L_i,~~ M_{2i}^{(1)}= M_2^{(1)} \cap L_i.
\end{equation*}

It is clear that the sets $M_{1i}^{(1)}$ and $M_{2i}^{(1)}$  are mapped to $M_{1i}^{(1)}$ and $M_{2i}^{(1)}$, respectively, setwise,  under each automorphism of $G$ preserving the labeling. Since  for any $0 \leq i \leq 2$, we have $|M_{1i}^{(1)}|\leq \lceil \sqrt[k]{k^{k -1}}\rceil -1$, so we can label all incident edges to each element of $M_{1i}^{(1)}$ with labels $\{1,2 \}$, such that for any two vertices of $M_{1i}^{(1)}$, say $x$ and $y$, there exists a label $c$, $1\leq c \leq 2$, such that the number of label $c$ for the incident edges to $x$ is different from the number of label $c$ for the incident edges to $y$. Hence, it can be deduce that each vertex of  $M_{1i}^{(1)}$ is fixed  under each automorphism of $G$ preserving the labeling, where $0\leq i \leq  2$. Thus every vertices of $M_1^{(1)}$ is fixed  under each automorphism of $G$ preserving the labeling. In sequel, we want to label the edges incident to vertices of $M_2^{(1)}$ such that $M_2^{(1)}$ is fixed  under each automorphism of $G$ preserving the labeling, pointwise. For this purpose, we partition the vertices of $M_{2i}^{(1)}$ to the sets $M_{{2i}_j}^{(1)}$ as follows ($1 \leq j \leq k -1$):
\begin{equation*}
M_{{2i}_j}^{(1)} =\left\{x\in M_{2i}^{(1)}~:~ |N(x) \cap  N^{(2)}(v)| = j \right\}.
\end{equation*}

Since the set $N^{(i)}(v)$, for any $i$, is mapped to itself, it can be concluded that $M_{{2i}_j}^{(1)}$  is mapped to itself under each automorphism of $G$ preserving the labeling, for any $i$ and $j$.  Let $M_{{2i}_j}^{(1)} = \{x_{j1}, x_{j2}, \ldots , x_{js_j}\}$. It is clear that $|M_{{2i}_j}^{(1)}| \leq |M_{2i}^{(1)}| \leq \lceil \sqrt[k]{k^{k -1}}\rceil-1$. Now we consider the two following cases for every $0 \leq i \leq 2$:

Case 1) Let  $1 \leq j < k -1$. Since   $|M_{{2i}_j}^{(1)}|\leq \lceil \sqrt[k]{k^{k -1}}\rceil-1$, so we can label all incident edges to each element of $M_{{2i}_j}^{(1)}$ with labels $\{1,2\}$, such that for any two vertices of $M_{{2i}_j}^{(1)}$, say $x$ and $y$, there exists a label $c$, $1\leq c \leq 2$, such that the number of label $c$ for the incident edges to $x$ is different from the number of label $c$ for the incident edges to $y$. Hence, it can be deduce that each vertex of  $M_{{2i}_j}^{(1)}$ is fixed  under each automorphism of $G$ preserving the labeling, where $1\leq j < k -1$.

Case 2) Let  $j = k -1$. Let $x_{jc}\in M_{{2i}_j}^{(1)}$, and $N(x_{jc}) \cap N^{(2)}(v) =\{x'_{jc1}, x'_{jc2}, \ldots , x'_{jcj}\}$. We assign to the  $j$-ary $(x_{jc}x'_{jc1}, \ldots , x_{jc}x'_{jcj})$ of edges, a $j$-ary of labels such that for every $x_{jc}$ and $x_{jc'}$, $1\leq c,c'\leq s_j$, there exists a label $l$ in their corresponding $j$-arys of labels with different number of label $l$ in their coordinates. For constructing $| M_{{2i}_j}^{(1)}|$ numbers of such $j$-arys we need, ${\rm min}\{r:~ {j+r-1 \choose r-1}\geq | M_{{2i}_j}^{(1)}| \}$ distinct labels. Since for any $j = k -1$, we have
\begin{equation*}
{\rm min}\left\{r:~ {j+r-1 \choose r-1}\geq | M_{{2i}_j}^{(1)}| \right\} \leq {\rm min}\left\{r:~ {j+r-1 \choose r-1}\geq \lceil \sqrt[k]{k^{k -1}}\rceil -1\right\} \leq 2,
\end{equation*}
so we need at most $2$ distinct labels from label set $\{1, 2\}$ for constructing such $j$-arys.   Hence, the vertices of $M_{{2i}_j}^{(1)}$, for any $j = k -1$, are fixed under each automorphism of $G$ preserving  the labeling.

Therefore, the vertices of $M_{2i}^{(1)}$ for any $0 \leq i \leq 2$, and so the vertices of $M_2^{(1)}$ are fixed under each automorphism of $G$ preserving  the labeling. Now, we can get that all vertices of $N^{(1)}(v)$ are fixed. If there exist unlabeled edges of $G$ with the two endpoints in $N^{(1)}(v)$, then we assign them an arbitrary label, say 1.

Step 2) Now we consider $N^{(2)}(v)$. We partition this set such that the vertices of $N^{(2)}(v)$ with the same neighbours in $M_2^{(1)}$, lie in a set. In other words, we can write $N^{(2)}(v) = \bigcup_i A_i$, such that $A_i$ contains that elements of $N^{(2)}(v)$ having the same neighbours in $M_2^{(1)}$, for any $i$. Since all vertices in  $M_2^{(1)}$  are fixed, so  the set $A_i$ is mapped to $A_i$ setwise, under each automorphism of $G$ preserving the labeling. Let $A_i = \{w_{i1}, \ldots , w_{it_i}\}$, and we have 
\begin{equation*}
N(w_{i1}) \cap M_2^{(1)} = \cdots = N(w_{it_i}) \cap M_2^{(1)} = \{v_{i1}, \ldots , v_{ip_i}\}.
\end{equation*}

We consider the two following cases:

Case 1) If for every $w_{ij}$ and $w_{ij'}$ in $A_i$, where $1\leq j,j' \leq t_i$, there exists a $c$, $1\leq c \leq p_i$, for which the label of edges $w_{ij}v_{ic}$ is different from label of edge $w_{ij'}v_{ic}$, then all vertices of $G$ in $A_i$ are fixed under each automorphism  of $G$ preserving the labeling.

Case 2) If there exist  $w_{ij}$ and $w_{ij'}$ in $A_i$, where $1\leq j,j' \leq t_i$, such that for every  $c$, $1\leq c \leq p_i$,  the label of edges $w_{ij}v_{ic}$ and $w_{ij'}v_{ic}$ are the same, then we can make a labeling  such that the  vertices  in $A_i$ have the same property   as Case 1, and so are fixed under each automorphism  of $G$ preserving the labeling, by using at least one of the following actions:
\begin{itemize}
\item By commutating the $j$-ary of labels assigned to the incident edges to $v_{ic}$ with an end point in $N^{(2)}(v)$, such that the vertices in $M_2^{(1)}$ are fixed under each automorphism of $G$ preserving the labeling,
\item By using a new $j$-ary of labels, with labels $\{1, 2\}$,  for incident edges to $v_{ic}$ with an end point in $N^{(2)}(v)$, such that the vertices in $M_2^{(1)}$ are fixed under each automorphism of $G$ preserving the labeling, 
\item  By labeling the unlabeled edges of $G$ with the two end points in $N^{(2)}(v)$ which are incident  to the vertices  in $A_i$, 
\item By labeling the unlabeled edges of $G$   which are incident  to the vertices  in $A_i$, and another their endpoint is $N^{(3)}(v)$,
\item By labeling the unlabeled edges of $G$ with the two end points in $N^{(3)}(v)$ for which the end points in $N^{(3)}(v)$ are adjacent  to some of  vertices  in $A_i$.
\end{itemize}

Using at least one of above actions, it can be seen that every two vertices $w_{ij}$ and $w_{ij'}$ in $A_i$ have the property as Case (1). Thus we conclude that all vertices in $A_i$, for any $i$, and so all vertices in $N^{(2)}(v)$, are fixed under each automorphism  of $G$ preserving the labeling.  If there exist unlabeled edges of $G$ with the two endpoints in $N^{(2)}(v)$, then we assign them an arbitrary label, say 1.

By following this method,  in the next step we  partition $N^{(3)}(v)$ exactly by the same method as partition of  $N^{(2)}(v)$ to the sets $A_i$s  in  Step 2, we can make a labeling such that $N^{(i)}(v)$ is fixed pointwise, under each automorphism  of $G$ preserving the labeling, for any $3 \leq i \leq {\rm diam}(G)$.\qed

 We end this  paper with proposing the following conjecture: 
 
 \begin{conjecture}
Let $G$ be a $k$-regular graph. If $k\geq 5$, then the distinguishing index of $G$ is at most two.
   \end{conjecture}
 
\end{document}